\documentclass[12pt, reqno]{amsart}
\title{2-Periodic Frieze Patterns}
\author{Alastair King \& Toril Palaniappan}
\date{30 June 2025}
\address{AK: Dept of Mathematical Sciences, University of Bath, Bath BA2 7AY, U.K.}
\email{a.d.king@bath.ac.uk} 
\address{TP: Dept of Mathematics, King's College London, Strand, London WC2R 2LS, U.K.}
\email{toril.palaniappan@kcl.ac.uk}
\usepackage[a4paper, hmargin=2.5cm,vmargin=2.5cm]{geometry}

\usepackage{amsfonts,amsmath,amssymb,amsthm,amsxtra}

\usepackage[usenames]{color}
\usepackage[colorlinks]{hyperref}
\usepackage{tikz}
  \usetikzlibrary{calc} 

\pgfmathsetmacro{\dotrad}{0.1}

\theoremstyle{definition}
\newtheorem{theorem}{Theorem}[section]

\newtheorem{lemma}[theorem]{Lemma}

\theoremstyle{definition}

\numberwithin{equation}{section}
\numberwithin{figure}{section}

\newcommand{\SL}{\operatorname{SL}}

\newcommand{\ZZ}{\mathbb{Z}}

\newcommand{\RR}{\mathbb{R}}
\newcommand{\RRp}{\RR^+}
\newcommand{\qint}[2]{\left[ #1 \right]_{#2}}
\newcommand{\qvar}{q}
\newcommand{\ev}{\mathrm{even}}
\newcommand{\od}{\mathrm{odd}}

\tikzset{%
friezeL/.style={xscale=0.8, yscale=0.7},
friezeM/.style={xscale=1.5, yscale=0.65}, %change this xscale to make it wider
friezeS/.style={xscale=0.7, yscale=0.6},
edge/.style={thick, blue, dashed},
}

\newcommand{\setfrzA}{
\pgfmathsetmacro{\frzLength}{4}
\pgfmathsetmacro{\frzhlfLength}{5}
\pgfmathsetmacro{\frzWidth}{7}
\pgfmathsetmacro{\frzTop}{1}
\pgfmathsetmacro{\frzBot}{-4}
\newcommand{\connect}{
  \foreach \x/\y in {a1/a2, a2/a3, a3/a4, a4/a5, a5/a6} \draw (\x)--(\y); }
}
\newcommand{\oddlayerA}[4]{
  \coordinate (atop) at (#1,\frzTop+\frzExtra);
  \coordinate (abot) at (#1,\frzBot-\frzExtra);  
\foreach \ver/\pos/\val in {a1/1/{#2}, a3/-1/{#3}, a5/-3/{#4}}
  \draw ({#1},\pos) node(\ver) {\val}; }
\newcommand{\evenlayerA}[4]{
  \coordinate (atop) at (#1,\frzTop+\frzExtra);
  \coordinate (abot) at (#1,\frzBot-\frzExtra);  
  \foreach \ver/\pos/\val in {a2/0/{#2}, a4/-2/{#3}, a6/-4/{#4}}
    \draw ({#1},\pos) node(\ver) {\val};
}

\newcommand{\setfrzDF}{% === D 4 ===
\pgfmathsetmacro{\frzLength}{7}
\pgfmathsetmacro{\frzhlfLength}{10}
\pgfmathsetmacro{\frzWidth}{7}
\pgfmathsetmacro{\frzTop}{1}
\pgfmathsetmacro{\frzBot}{-1}
\newcommand{\connect}{
\foreach \x/\y in {a1/a4, a2/a4, a3/a4} \draw (\x)--(\y); }
}
\newcommand{\oddlayerDF}[4]{
  \coordinate (atop) at (#1,\frzTop+\frzExtra);
  \coordinate (abot) at (#1,\frzBot-\frzExtra);  
\foreach \ver/\pos/\val in {a1/1/{#2}, a2/0/{#3}, a3/-1/{#4}}
  \draw ({#1},\pos) node(\ver) {\val}; 
}
\newcommand{\evenlayerDF}[2]{
  \coordinate (atop) at (#1,\frzTop+\frzExtra);
  \coordinate (abot) at (#1,\frzBot-\frzExtra);  
  \foreach \ver/\pos/\val in {a4/0/{#2}}
    \draw ({#1},\pos) node(\ver) {\val}; 
}
\newcommand{\setfrzDE}{
\pgfmathsetmacro{\frzLength}{7}
\pgfmathsetmacro{\frzhlfLength}{10}
\pgfmathsetmacro{\frzWidth}{7}
\pgfmathsetmacro{\frzTop}{5}
\pgfmathsetmacro{\frzBot}{-1}
\newcommand{\connect}{
  \foreach \x/\y in {a1/a2, a3/a2, a3/a4, a5/a4, a5/a6, a7/a6, a9/a6} \draw (\x)--(\y); }
}

\newcommand{\setfrzDO}{
\pgfmathsetmacro{\frzTop}{6}
\pgfmathsetmacro{\frzBot}{-1}
\newcommand{\connect}{
  \foreach \x/\y in {a0/a1, a1/a2, a3/a2, a3/a4, a5/a4, a5/a6, a7/a6, a9/a6} \draw (\x)--(\y); }
}

\pgfmathsetmacro{\frzExtra}{0.7}
\newcommand{\perEdge}[1]{
   \draw [edge] (#1,\frzTop+\frzExtra)--(#1,\frzBot-\frzExtra); }

\begin{document}  \maketitle
% ================================================================
\begin{abstract}
We classify 2-periodic mesh friezes of finite type $A$, $D$ or $E$
with positive real entries.
There are families with 0,1, or 2 parameters, depending on type.
\end{abstract}
% ================================================================
\section{Introduction}
% ================================================================

Frieze patterns were introduced and studied by Coxeter and Conway \cite{Cox,CC}.
Such frieze patterns would now be considered to be of type $A$
and can naturally be generalised to other quiver types, 
particularly other Dynkin types $D$ and $E$,
where a similar periodicity remains
(see e.g.~\cite[\S2]{MG}).

In the context of cluster algebra, it is also very natural to consider friezes with values in $\RRp$,
the set of positive real numbers.
In particular, in Dynkin types, the set of all $\RRp$-friezes is the set of positive points of a 
cluster variety of finite type.
This variety carries a lot of interesting structure (see e.g.~\cite{GM}).
In particular, it carries the action of a finite order automorphism $\tau$
coming from translating the frieze one step to the left
(also recognised as the DT-transformation \cite{SGL}).
It turns out that the $\tau$-action has a unique fixed point (see e.g.~\cite[\S3]{SGL} for more detail),
that is, there is a unique $\RRp$-valued constant (or 1-periodic) frieze.
For example, in type $D_4$, the unique constant frieze is 
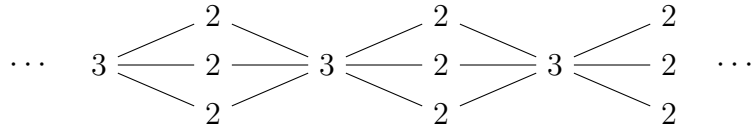
\[ %\begin{figure}[ht]
\begin{tikzpicture}[friezeM] \setfrzDF % == D 4  ==
%\perEdge{0.5}; \perEdge{4.5};
\evenlayerDF{0}{3}
\oddlayerDF{1}{2}{2}{2} \connect
\evenlayerDF{2}{3} \connect 
\oddlayerDF{3}{2}{2}{2} \connect 
\evenlayerDF{4}{3} \connect 
\oddlayerDF{5}{2}{2}{2} \connect
\draw (-0.6,0) node {$\cdots$};
\draw (5.6,0) node {$\cdots$};
\end{tikzpicture}
\] %\end{figure}
However, in general, the constant frieze is rarely integral.

The objective of this article is to classify 2-periodic $\RRp$-friezes. 
We will see that the behaviour in this case is type dependent.
As an example, consider a 2-periodic $\RRp$-frieze of type $D_4$, that is, the following pattern.
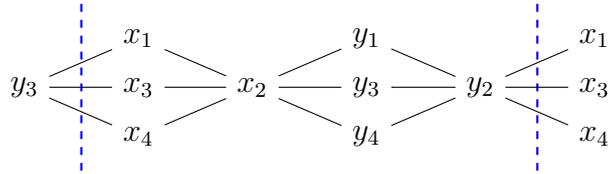
\[ %\begin{figure}[ht]
\begin{tikzpicture}[friezeM] \setfrzDF % == D 4  ==
\perEdge{0.5}; \perEdge{4.5};
\evenlayerDF{0}{$y_3$}
\oddlayerDF{1}{$x_1$}{$x_3$}{$x_4$} \connect
\evenlayerDF{2}{$x_2$} \connect 
\oddlayerDF{3}{$y_1$}{$y_3$}{$y_4$} \connect 
\evenlayerDF{4}{$y_2$} \connect 
\oddlayerDF{5}{$x_1$}{$x_3$}{$x_4$} \connect
\end{tikzpicture}
\] %\end{figure}
in which the variables $x_i$ and $y_i$ take values in $\RRp$.
The frieze equations are
\begin{align*}
x_1y_1 = x_3y_3 = x_4y_4 &= x_2 + 1 = y_2 + 1 \\
%x_4y_4 &= x_3 + 1 = y_3 + 1 \\
%x_2y_2 &= x_3 + 1 = y_3 + 1 \\
%x_4y_4 &= x_3+ 1 = y_3 + 1 \\
x_2y_2 &= x_1x_3x_4 + 1 = y_1y_3y_4 + 1
\end{align*}
These immediately give $x_2=y_2 =c$, say, and then $(c^2-1)^2= (x_1x_2x_4)(y_1y_2y_4)=(c+1)^3$,
that is, $(c-3)c(c+1)^2=0$.
Hence the only positive solution is $c=3$
and the frieze equations simplify to
 \begin{align*}
x_1y_1 = x_3y_3 = x_4y_4 = 4
\qquad  x_1x_3x_4 = y_1y_3y_4 \,(=8)
\end{align*}
so that two free parameters, e.g.~$x_1,x_3\in\RRp$, determine a 2-periodic frieze of type $D_4$.

It turns out that this example exhibits the behaviour for a general 2-periodic $\RRp$-frieze of finite type.
Namely, some $\tau$-orbits are constant and on the other $\tau$-orbits, the two values $x,y$ satisfy
$xy=\lambda$, where $\sqrt{\lambda}$ is the value on that orbit in the constant frieze.

In general,  we will prove the following.

\begin{theorem}\label{thm:main}
The classification of 2-periodic $\RRp$-friezes of types $A$, $D$, $E$ is as follows.
\begin{enumerate}
\item 
All 2-periodic friezes of types $A_\ev$, $E_6$ and $E_8$ are constant.
\item
There are one-parameter families of 2-periodic friezes of types $A_\od$, $D_\od$ and $E_7$.
\item
There are two-parameter families of 2-periodic friezes of types $D_\ev$. 
\end{enumerate}
\end{theorem}

The proofs and more precise statements are given in Sections 3--5.
Note that, for types $A$ and $E$, the arguments will exploit the fact 
that there is a unique constant frieze of each type.
On the other hand, for type $D$, the classification is self-contained and thus reproves this uniqueness.

Observe also that Theorem~\ref{thm:main} is consistent with \cite[Thm~3.12]{SGL}, 
in that the number of parameters is equal to the multiplicity of the $-1$-eigenvalue of 
the derivative/linearisation of $\tau$ at the fixed point,
that is, the number of Coxeter exponents $e$ for which $2e=h$, the Coxeter number,
or rather $e+1=(h+2)/2$.

% ===============================================
\section{Constant Friezes}\label{sec:const}
% ===============================================
\pgfmathsetmacro{\m}{2*cos(180/4)} % {sqrt(2)}
\newcommand{\msym}{\sqrt{2}}
\pgfmathsetmacro{\p}{2*cos(180/5)} % {\phi= (1+sqrt(5))/2}
\newcommand{\psym}{\phi}
\pgfmathsetmacro{\q}{2*cos(180/6)} % {sqrt(3)}
\newcommand{\qsym}{\sqrt{3}}
\pgfmathsetmacro{\r}{2*cos(180/7)}
\pgfmathsetmacro{\rr}{\r^2-1}
\newcommand{\rsym}{\alpha}
\newcommand{\rrsym}{\beta}
\newcommand{\eng}{\ell}

Constant friezes of type $A_n$ were known by Conway--Coxeter \cite[Problems 13-16]{CC}
to be given by the lengths of the diagonals of a regular $\eng$-gon of side 1,
where $\eng=n+3$.

\begin{equation}\label{eq:h-gons}
\pgfmathsetmacro{\voffset}{-1.2}
% ===========
\begin{tikzpicture} [scale=1.3, baseline=(bb.base)] \coordinate (bb) at (0,0);
  \pgfmathsetmacro{\N}{4}
  \pgfmathsetmacro{\ang}{360/\N}
  \pgfmathsetmacro{\circrad}{0.5*cosec(0.5*\ang)}
  \pgfmathsetmacro{\stang}{180}
\foreach \n in {0,...,\N}
  \coordinate (v\n) at (\stang-\n*\ang:\circrad);
\foreach \n in {1,...,\N}
  \pgfmathsetmacro{\nm}{\n-1}
  \draw (v\nm) to (v\n);
  \pgfmathsetmacro{\M}{\N-1}
\foreach \n in {1,...,\M}
  \draw [blue, thick] (v0) to (v\n);
%\foreach \n in {3,...,\M}
%  \draw [blue, thick, dotted] (v1) to (v\n);
\foreach \n/\lab in {2/\msym}
  \draw ($0.5*(v0)+0.5*(v\n)$) node [above] {\small $\lab$};
\draw (0,\voffset) node {\small $\msym \approx 1.414$};
 \end{tikzpicture} 
% ===========
 \qquad
% ===========
 \begin{tikzpicture} [scale=1.3, baseline=(bb.base)] \coordinate (bb) at (0,0);
  \pgfmathsetmacro{\N}{5}
  \pgfmathsetmacro{\ang}{360/\N}
  \pgfmathsetmacro{\circrad}{0.5*cosec(0.5*\ang)}
  \pgfmathsetmacro{\stang}{180}
\foreach \n in {0,...,\N}
  \coordinate (v\n) at (\stang-\n*\ang:\circrad);
\foreach \n in {1,...,\N}
  \pgfmathsetmacro{\nm}{\n-1}
  \draw (v\nm) to (v\n);
  \pgfmathsetmacro{\M}{\N-1}
\foreach \n in {1,...,\M}
  \draw [blue, thick] (v0) to (v\n);
%\foreach \n in {3,...,\M}
%  \draw [blue, thick, dotted] (v1) to (v\n);
\foreach \n/\lab in {2/\psym, 3/\psym}
  \draw ($0.5*(v0)+0.5*(v\n)$) node [above] {\small $\lab$};
\draw (0,\voffset) node {\small $\psym \approx 1.618$};
 \end{tikzpicture} 
% ===========
 \qquad
% ===========
  \begin{tikzpicture} [scale=1.3, baseline=(bb.base)] \coordinate (bb) at (0,0);
  \pgfmathsetmacro{\N}{6}
  \pgfmathsetmacro{\circrad}{1}
  \pgfmathsetmacro{\ang}{360/\N}
  \pgfmathsetmacro{\stang}{180}
\foreach \n in {0,...,\N}
  \coordinate (v\n) at (\stang-\n*\ang:\circrad);
\foreach \n in {1,...,\N}
  \pgfmathsetmacro{\m}{\n-1}
  \draw (v\m) to (v\n);
  \pgfmathsetmacro{\M}{\N-1}
\foreach \n in {1,...,\M}
  \draw [blue, thick] (v0) to (v\n);
%\foreach \n in {3,...,\M}
%  \draw [blue, thick, dotted] (v1) to (v\n);
\foreach \n/\lab in {2/\qsym, 3/2, 4/\qsym}
  \draw ($0.5*(v0)+0.5*(v\n)$) node [above] {\small $\lab$};
\draw (0,\voffset) node {\small $\qsym \approx 1.732$};
 \end{tikzpicture} 
% ===========
 \quad
% ===========
  \begin{tikzpicture} [scale=1.3, baseline=(bb.base)] \coordinate (bb) at (0,0);
  \pgfmathsetmacro{\N}{7}
  \pgfmathsetmacro{\circrad}{1}
  \pgfmathsetmacro{\ang}{360/\N}
  \pgfmathsetmacro{\stang}{180}
\foreach \n in {0,...,\N}
  \coordinate (v\n) at (\stang-\n*\ang:\circrad);
\foreach \n in {1,...,\N}
  \pgfmathsetmacro{\m}{\n-1}
  \draw (v\m) to (v\n);
  \pgfmathsetmacro{\M}{\N-1}
\foreach \n in {1,...,\M}
  \draw [blue, thick] (v0) to (v\n);
%\foreach \n in {3,...,\M}
%  \draw [blue, thick, dotted] (v1) to (v\n);
\foreach \n/\lab in {2/\rsym, 3/\rrsym, 4/\rrsym, 5/\rsym}
  \draw ($0.5*(v0)+0.5*(v\n)$) node [above] {\small $\lab$};
\draw (0,\voffset) node {\small $\rsym \approx 1.802,\; \rrsym \approx 2.247$};
%\draw (0.1,1.4*\voffset) node {\small $\rrsym \approx 2.247$};
 \end{tikzpicture} 
\end{equation}
% ===========

The length of the $k$-th diagonal of an $\eng$-gon of side 1 is explicitly given by
\[
  \qint{k}{\eng} = \frac{\sin(k\pi /\eng)}{\sin(\pi /\eng)} = \frac{\qvar^k-\qvar^{-k}}{\qvar-\qvar^{-1}}
  \quad\text{for $\qvar=e^{i\pi/\eng}$,}
\]
that is, quantum integers (or $\SL(2)$-characters) at a $2\eng$-th root of unity.
Note that the `1st diagonal' is the side, so the first proper diagonal is the 2nd,
of length $\qint{2}{\eng}=2\cos (\pi/\eng)$.
Alternatively, the quantum integers are given by $\qint{k+1}{\eng} = p_{k}(\qint{2}{\eng})$,
where the Chebyshev polynomials $p_{n}(x)$ are defined recursively by
$p_{m+1}=x p_m -p_{m-1}$, starting from $p_0=1$ and $p_{1}=x$. % $p_1=x$.

We get the most general infinite constant frieze this way, 
because Chebyshev polynomials also satisfy
\begin{equation}\label{eq:cluster-rec}
  p_m^2 =  p_{m+1} p_{m-1} +1,
\end{equation}
which gives an alternative defining recurrence.
Setting $x= \qint{2}{n+3}$ %, for $h=n+3$, 
gives $p_{n+1}(x)=1$ and $p_{n+2}(x)=0$,
while $p_m(x)>0$ for $1\leq m\leq n$,
giving a constant frieze of type $A_n$.

Constant friezes of type $D_n$ are easy to write down, with $2,3,\ldots,n-1$ on the long leg 
and $\sqrt{n}$ at the ends of the short legs.
For example:

\begin{equation}\label{eq:constD}
\tikzset{edge/.style={thick,blue},
vert/.style={fill=blue}}
% ==== D4 ====
\begin{tikzpicture} [scale=0.75, baseline=(bb.base)] 
\coordinate (bb) at (0,0);
\pgfmathsetmacro{\NN}{4}
\pgfmathsetmacro{\N}{\NN-2}
\coordinate (X) at (0,0); 
\foreach \n in {1} \coordinate (Va\n) at (-\n,0); 
\foreach \n in {1} \coordinate (Vb\n) at (0,-\n);
\foreach \n in {1,...,\N} \coordinate (Vc\n) at (\N-\n,0);
\draw [edge] (Va1)--(Vc1) (Vb1)--(X);
\foreach \n/\val in {1/2} 
  \draw [vert] (Va\n) circle (\dotrad) node [above] {\small $\val$}; 
\foreach \n/\val in {1/2} 
  \draw [vert] (Vb\n) circle (\dotrad) node [left] {\small $\val$};
\foreach \n/\val in {1/2, 2/3} 
  \draw [vert] (Vc\n) circle (\dotrad) node [above] {\small $\val$};
\end{tikzpicture}
\qquad
% ==== D5 ====
\begin{tikzpicture} [scale=0.75, baseline=(bb.base)] 
\coordinate (bb) at (0,0);
\pgfmathsetmacro{\NN}{5}
\pgfmathsetmacro{\N}{\NN-2}
\newcommand{\sq}{\sqrt{\NN}}
\coordinate (X) at (0,0); 
\foreach \n in {1} \coordinate (Va\n) at (-\n,0); 
\foreach \n in {1} \coordinate (Vb\n) at (0,-\n);
\foreach \n in {1,...,\N} \coordinate (Vc\n) at (\N-\n,0);
\draw [edge] (Va1)--(Vc1) (Vb1)--(X);
\foreach \n/\val in {1/\sq} 
  \draw [vert] (Va\n) circle (\dotrad) node [left] {\small $\val$}; 
\foreach \n/\val in {1/\sq} 
  \draw [vert] (Vb\n) circle (\dotrad) node [left] {\small $\val$};
\foreach \n/\val in {1/2, 2/3, 3/4} 
  \draw [vert] (Vc\n) circle (\dotrad) node [above] {\small $\val$};
\end{tikzpicture}
\qquad
% ==== D6 ====
\begin{tikzpicture} [scale=0.75, baseline=(bb.base)] 
\coordinate (bb) at (0,0);
\pgfmathsetmacro{\NN}{6}
\pgfmathsetmacro{\N}{\NN-2}
\newcommand{\sq}{\sqrt{\NN}}
%\coordinate (X) at (0,0); 
\foreach \n in {1} \coordinate (Va\n) at (-\n,0); 
\foreach \n in {1} \coordinate (Vb\n) at (0,-\n);
\foreach \n in {1,...,\N} \coordinate (Vc\n) at (\N-\n,0);
\draw [edge] (Va1)--(Vc1) (Vb1)--(X);
\foreach \n/\val in {1/\sq} 
  \draw [vert] (Va\n) circle (\dotrad) node [left] {\small $\val$}; 
\foreach \n/\val in {1/\sq} 
  \draw [vert] (Vb\n) circle (\dotrad) node [left] {\small $\val$};
\foreach \n/\val in {1/2, 2/3, 3/4, 4/5} 
  \draw [vert] (Vc\n) circle (\dotrad) node [above] {\small $\val$};
\end{tikzpicture}
\end{equation}

\goodbreak
Constant friezes of type $E$ can be found with some work, but the solutions can easily be verified
and are as follows.

\begin{equation}\label{eq:E6const}
% ===== E 6 algebraic =====
\begin{tikzpicture} [scale=0.5, baseline=(bb.base),
  edge/.style={thick,blue},
  vert/.style={fill=blue}]
\coordinate (bb) at (0,0);
\draw (-7,0) node {$E_6:$};
\pgfmathsetmacro{\NN}{6}
\pgfmathsetmacro{\N}{\NN-3}
\newcommand{\sBa}{\beta}
\newcommand{\sBb}{\alpha+\beta}
\newcommand{\sAa}{1+\alpha}
\newcommand{\sCa}{\beta}
\newcommand{\sCb}{\alpha+\beta}
\newcommand{\sXX}{1+2\alpha+\beta}
\coordinate (X) at (0,0); 
\foreach \n in {1,2} \coordinate (Va\n) at (-2+\n,0); 
\foreach \n in {1,2,3} \coordinate (Vb\n) at (0,3-\n);
\foreach \n in {1,...,\N} \coordinate (Vc\n) at (0,-\N+\n);
%\draw (0,3.3) node {\small $E_\NN$};
\draw [edge] (Vb1)--(Vc1) (Va1)--(X);
\draw [vert] (X) circle (\dotrad) 
  node [right] {\small $x=\sXX$};
\foreach \n/\val in {1/\sAa} 
  \draw [vert] (Va\n) circle (\dotrad) node [left] {\small $\val=a_{\n}$}; 
\foreach \n/\val in {1/\sBa, 2/\sBb} 
  \draw [vert] (Vb\n) circle (\dotrad) node [right] {\small $b_{\n} =\val$};
\foreach \n/\val in {1/\sCa, 2/\sCb} 
  \draw [vert] (Vc\n) circle (\dotrad) node [right] {\small $c_{\n} =\val$};
\end{tikzpicture}
\end{equation}
where $\alpha,\beta$ are the diagonals of a heptagon, which satisfy
\begin{equation}\label{eq:alpbet-too}
\begin{aligned}
  \alpha^2 &= 1+\beta, \qquad
  \beta^2 &=  1+\alpha+ \beta, \qquad
  \alpha\beta &= \alpha+ \beta.
\end{aligned}
\end{equation}

%For $E_7$, the constant frieze is:
\begin{equation}\label{eq:E7const}
% ===== E 7 algebraic =====
\begin{tikzpicture} [scale=0.5, baseline=(bb.base),
  edge/.style={thick,blue},
  vert/.style={fill=blue}]
\coordinate (bb) at (0,0);
\draw (-9,0) node {$E_7:$};
\pgfmathsetmacro{\NN}{7}
\pgfmathsetmacro{\N}{\NN-3}
\newcommand{\sBa}{1+\phi} 
\newcommand{\sBb}{1+3\phi} 
\newcommand{\sAa}{\sqrt{2}(1+\phi)} 
\newcommand{\sCa}{\sqrt{2}\phi} 
\newcommand{\sCb}{1+2\phi} 
\newcommand{\sCc}{2\sqrt{2}(1+\phi)}) 
\newcommand{\sXX}{3(1+2\phi)} 
\coordinate (X) at (0,0); 
\foreach \n in {1,2} \coordinate (Va\n) at (-2+\n,0); 
\foreach \n in {1,2,3} \coordinate (Vb\n) at (0,3-\n);
\foreach \n in {1,...,\N} \coordinate (Vc\n) at (0,-\N+\n);
%\draw (0,3.5) node {\small $E_\NN$};
\draw [edge] (Vb1)--(Vc1) (Va1)--(X);
% ===
\draw [vert] (X) circle (\dotrad) 
  node [right] {\small $x=\sXX$};
\foreach \n/\val in {1/\sAa}
  \draw [vert] (Va\n) circle (\dotrad) node [left] {\small $\val=a_{\n}$}; 
\foreach \n/\val in {1/\sBa, 2/\sBb} 
  \draw [vert] (Vb\n) circle (\dotrad) node [right] {\small $b_{\n} =\val$};
\foreach \n/\val in {1/\sCa, 2/\sCb, 3/\sCc} 
  \draw [vert] (Vc\n) circle (\dotrad) node [right] {\small $c_{\n} =\val$};
\end{tikzpicture}
\end{equation}
where $\phi$ is the golden ratio, i.e.~the diagonal of a pentagon,
which satisfies $\phi^2=\phi+1$.

%For $E_8$, the constant frieze is:
\begin{equation}\label{eq:E8const}
% ===== E 8 algebraic =====
\begin{tikzpicture} [scale=0.5, baseline=(bb.base),
  edge/.style={thick,blue},
  vert/.style={fill=blue}]
\coordinate (bb) at (0,0);
\draw (-8,0) node {$E_8:$};
\pgfmathsetmacro{\NN}{8}
\pgfmathsetmacro{\N}{\NN-3}
\newcommand{\sBa}{2+\sqrt{2}} 
\newcommand{\sBb}{5+4\sqrt{2}}
\newcommand{\sAa}{3+2\sqrt{2}} 
\newcommand{\sCa}{1+\sqrt{2}} 
\newcommand{\sCb}{2+2\sqrt{2}} 
\newcommand{\sCc}{5+3\sqrt{2}}  
\newcommand{\sCd}{9+6\sqrt{2}} 
\newcommand{\sXX}{16+12\sqrt{2}}  
\coordinate (X) at (0,0); 
\foreach \n in {1,2} \coordinate (Va\n) at (-2+\n,0); 
\foreach \n in {1,2,3} \coordinate (Vb\n) at (0,3-\n);
\foreach \n in {1,...,\N} \coordinate (Vc\n) at (0,-\N+\n);
%\draw (0,3.5) node {\small $E_\NN$};
\draw [edge] (Vb1)--(Vc1) (Va1)--(X);
% ===
\draw [vert] (X) circle (\dotrad) 
 node [right] {\small $x=\sXX$};
\foreach \n/\val in {1/\sAa} 
  \draw [vert] (Va\n) circle (\dotrad) node [left] {\small $\val=a_{\n}$}; 
\foreach \n/\val in {1/\sBa, 2/\sBb} 
  \draw [vert] (Vb\n) circle (\dotrad) node [right] {\small $b_{\n} =\val$};
\foreach \n/\val in {1/\sCa, 2/\sCb, 3/\sCc, 4/\sCd} 
  \draw [vert] (Vc\n) circle (\dotrad) node [right] {\small $c_{\n} =\val$};
\end{tikzpicture}
\end{equation}

Constant $\RRp$-friezes are special cases of $Q$-systems 
and it is known that these admit unique solutions in $\RRp$: 
see \cite[\S14]{KNS}, \cite[\S1]{NK} and also \cite[\S3]{SGL}.

% ===============================================
\section{2-Periodic Friezes of Type $A$}
% ===============================================

Similar to constant friezes of type $A$, we can describe 2-periodic friezes in terms of 
\emph{two-coloured quantum integers} (cf.~\cite[\S2]{Haz}).
These are two sequences of two-variable polynomials $p^x_n,\, p^y_n\in\ZZ[x,y]$, 
defined recursively by
\begin{equation}
  \begin{aligned}
    p^x_{n+1} &= x\, p^y_n - p^x_{n-1}  \\
    p^y_{n+1} &= y\, p^x_n - p^y_{n-1}
  \end{aligned}
\end{equation}
starting from $p_{0}^x =  p_{0}^y = 1$, while $p_1^x= x$ and $p_1^y= y$.
We can also set $p_{-1}^x =  p_{-1}^y = 0$.
The sequences continue
\begin{equation}\label{eq:first-p}
  \begin{aligned}
%    p^x_{1} &= x  &\qquad p^y_{1} &=y \\
    p^x_{2} &= xy-1  &\qquad p^y_{2} &=xy-1 \\
    p^x_{3} &= x(xy-2)  &\qquad p^y_{3} &=y(xy-2) \\
    p^x_{4} &= (xy)^2-3xy+1  &\qquad p^y_{4} &=(xy)^2-3xy+1
  \end{aligned}
\end{equation}
and, in general, setting $\lambda=xy$, we have
\begin{equation}\label{eq:gen-p}
  \begin{aligned}
   p^x_{2n-1} & =x\, q_{2n-1}(\lambda) &\qquad \quad p^y_{2n-1} &= y\, q_{2n-1}(\lambda) \\
    p^x_{2n}  &= q_{2n}(\lambda) &\qquad p^y_{2n} &= q_{2n}(\lambda)
   \end{aligned}
\end{equation}
for certain one-variable polynomials $q_m(\lambda)$.

As in the case of ordinary quantum integers/Chebyshev polynomials,
we have the following alternative recurrence.

\begin{lemma}
For $n\geq 1$, we have 
\begin{equation}\label{eq:frz-for-p}
   p^x_{n}p^y_{n-2} = p^x_{n-1}p^y_{n-1} - 1 = p^y_{n}p^x_{n-2}. 
\end{equation}
\end{lemma}

\begin{proof}
This follows by induction, with the case $n=1$ being immediate from the initial conditions.
Assuming \eqref{eq:frz-for-p} holds for $n$, we have
\begin{align*}
   p^x_{n+1}p^y_{n-1}
  & =  (x p^y_n - p^x_{n-1}) p^y_{n-1} \\
  & =  x p^y_n p^y_{n-1} - p^y_{n}p^x_{n-2} - 1\\
  & =  p^y_n (x p^y_{n-1} - p^x_{n-2}) - 1\\
  & =  p^y_n  p^x_n - 1
\end{align*}
so that the left-hand equality holds for $n+1$, 
while the right-hand equality follows by exchanging the roles of $x$ and $y$.
\end{proof}

Thus these two sequences form a general infinite 2-periodic frieze when embedded diagonally as follows:
\begin{equation}\label{eq:typeAfrz}
\begin{tikzpicture}[friezeM, baseline=(bb.base)] \setfrzA % ==  A  ==
\perEdge{0.5};   \perEdge{4.5}; 
\evenlayerA{0}{$p_2^y$}{$p_4^x$}{$\vdots$}
\oddlayerA{1}{$p_1^x$}{$p_3^y$}{$p_5^x$} \connect
\evenlayerA{2}{$p_2^x$}{$p_4^y$}{$\vdots$} \connect 
\oddlayerA{3}{$p_1^y$}{$p_3^x$}{$p_5^y$} \connect 
\evenlayerA{4}{$p_2^y$}{$p_4^x$}{$\vdots$} \connect 
\oddlayerA{5}{$p_1^x$}{$p_3^y$}{$p_5^x$} \connect
\coordinate (bb) at (0,-1);
\end{tikzpicture}
\end{equation}
In particular, by \eqref{eq:gen-p}, the even layers are constant.

To obtain a 2-periodic $\RRp$-frieze of width $w$ (or type $A_w$), 
we need $x, y\in\RRp$ with $p^x_{w+1}=p^y_{w+1}=1$, 
while $p^x_{n}>0$ and $p^y_{n}>0$, for $1\leq n\leq w$.
If $w$ is even, then \eqref{eq:gen-p} imples that $x=y$ 
and so we have a constant frieze, that is,
\[
  p^x_{n}=p^y_{n} = \qint{n+1}{h} \quad\text{for $h=w+3$}\
\]
If $w$ is odd, then we need to solve $q_{w+1}(\lambda)=1$ for $\lambda>0$ 
and $x,y$ are then just constrained by $xy=\lambda$.
Since setting $x=y=\sqrt{\lambda}$ gives a constant frieze,
we see that $\lambda$ must take the value it does for the constant frieze of the same width,
namely $\sqrt{\lambda}=\qint{2}{h}$,
where again $h=w+3$.
Thus we have a one-parameter family of 2-periodic friezes as in \eqref{eq:typeAfrz}, with
\begin{equation}\label{eq:last-p}
  \begin{aligned}
    p^x_{n} &= \qint{n+1}{h}   &\quad p^y_{n} &= \qint{n+1}{h} &\quad\text{for $n$ even,}\\
    p^x_{n} &= x \qint{n+1}{h} / \qint{2}{h}  &\quad
    p^y_{n} &= y \qint{n+1}{h} / \qint{2}{h} &\quad \text{for $n$ odd.}
  \end{aligned}
\end{equation}

% ===============================================
\section{2-Periodic Friezes of Type $D$}
% ===============================================

Now consider a 2-periodic frieze of type $D_n$ with $n$ even, say $n=2m+2$.
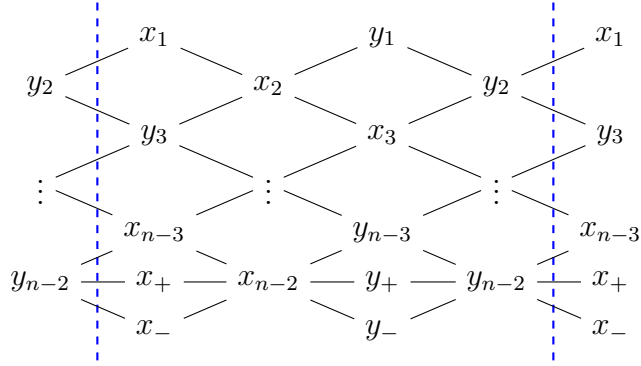
\[ %\begin{figure}[ht]
\begin{tikzpicture}[friezeM] \setfrzDE % ==  D Even  ==
\newcommand{\oddlayerDE}[6]{
  \coordinate (atop) at (#1,\frzTop+\frzExtra);
  \coordinate (abot) at (#1,\frzBot-\frzExtra);  
\foreach \ver/\pos/\val in {a1/5/{#2}, a3/3/{#3}, a5/1/{#4}, a7/0/{#5}, a9/-1/{#6}}
  \draw ({#1},\pos) node(\ver) {\val}; }
\newcommand{\evenlayerDE}[4]{
  \coordinate (atop) at (#1,\frzTop+\frzExtra);
  \coordinate (abot) at (#1,\frzBot-\frzExtra);  
  \foreach \ver/\pos/\val in {a2/4/{#2}, a4/2/{#3}, a6/0/{#4}}
    \draw ({#1},\pos) node(\ver) {\val}; }
\perEdge{0.5}; \perEdge{4.5};
\evenlayerDE{0}{$y_2$}{\vdots}{$y_{n-2}$}
\oddlayerDE{1}{$x_1$}{$y_3$}{$x_{n-3}$}{$x_{+}$}{$x_{-}$} \connect
\evenlayerDE{2}{$x_2$}{\vdots}{$x_{n-2}$} \connect 
\oddlayerDE{3}{$y_1$}{$x_3$}{$y_{n-3}$}{$y_{+}$}{$y_{-}$} \connect 
\evenlayerDE{4}{$y_2$}{\vdots}{$y_{n-2}$} \connect 
\oddlayerDE{5}{$x_1$}{$y_3$}{$x_{n-3}$}{$x_{+}$}{$x_{-}$} \connect
\end{tikzpicture}
\] %\end{figure}
The frieze equations are
\begin{equation}\label{eq:DEfrzeq}
\begin{aligned}
x_1y_1 &= x_2+ 1 = y_2 + 1 \\
x_2y_2 &= x_1y_3 + 1 = y_1x_3 + 1 \\
x_3y_3 &= x_2y_4 + 1 = y_2x_4 + 1 \\
& \;\;\vdots \\ % \nonumber
x_{n-3}y_{n-3} &= x_{n-2}y_{n-4}+1 = y_{n-2}x_{n-4}+1 \\
x_{n-2}y_{n-2} &= x_{+}x_{-}x_{n-3} + 1 = y_{+}y_{-}y_{n-3} + 1 \\ 
x_{+}y_{+}  = x_{-}y_{-} &= x_{n-2} + 1 = y_{n-2} + 1
\end{aligned}
\end{equation}
These equations imply $x_{2k}= y_{2k}(=c_{2k})$, for $k = 1,\dots,m$.
We then get the following equations for the $c_{2k}$, including the effective initial value $c_0=1$,
\begin{equation}\label{eq:c2k}
\begin{aligned}
(c_{2m}^2-1)^2 &= (x_{+}x_{-}x_{2m-1})(y_{+}y_{-}y_{2m-1}) = (c_{2m}+1)^2(c_{2m}c_{2m-2}+1)\\
(c_{2k}^2-1)^2 &= (x_{2k+1}y_{2k-1})(y_{2k+1}x_{2k-1}) = (c_{2k+2}c_{2k}+1)(c_{2k}c_{2k-2}+1)
\end{aligned}
\end{equation}
for $k=1,\dots,m-1$.
From these, starting from $k=m$, we get $c_{2k-2}=c_{2k}-2$.
Then, starting from $c_0=1$, we conclude that $c_{2k}=2k+1$.
Putting this back into \eqref{eq:DEfrzeq} gives
\begin{equation}\label{eq:xyodd}
\begin{aligned}
  x_{2k-1}y_{2k-1} &= (2k-1)(2k+1)+1 =4k^2, \\
   x_{2k-1}y_{2k+1} = x_{2k+1}y_{2k-1} &= (2k+1)^2-1 = 4k(k+1)
\end{aligned}
\end{equation}
and hence, inductively, $x_{2k-1}=kx_1$, $y_{2k-1}=ky_1$.
Thus we have reduced the frieze equations \eqref{eq:DEfrzeq} to 
\begin{equation}\label{eq:DEfrz-final}
 x_1y_1=4,\quad x_{+}y_{+}  = x_{-}y_{-} = n,\quad  x_{1}x_{+}x_{-} = y_{1}y_{+}y_{-} (=2n),
\end{equation}
where the equation in parentheses follows because
$(x_{1}x_{+}x_{-})(y_{1}y_{+}y_{-})=4n^2$ from the first three equations.
%This shows that we have two free parameters, just as in the $D_4$ case.
As shown above, the other values in the frieze are
\begin{equation}\label{eq:DEfrz-other}
 x_{2k}=y_{2k}=2k+1,\quad x_{2k-1}=kx_1,\quad y_{2k-1}=ky_1.
\end{equation}
Thus there are two free parameters, e.g.~given by $\{(x_+,y_+)\in(\RRp)^2 : x_+y_+=n\}$
and $\{(x_-,y_-)\in(\RRp)^2 : x_-y_-=n\}$,
Setting $x=2n/(x_+x_-)$ and $y=2n/(y_+y_-)$
 the frieze is
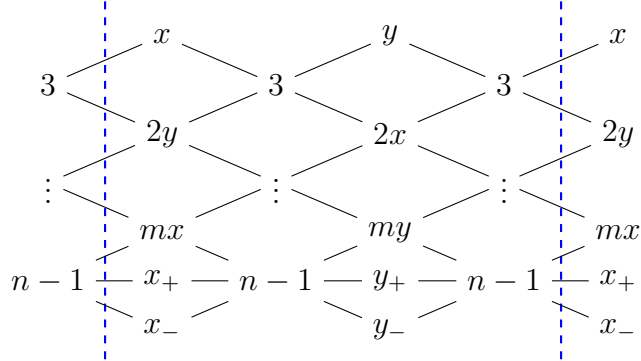
\[ %\begin{figure}[ht]
\begin{tikzpicture}[friezeM] \setfrzDE % ==  D Even  ==
\newcommand{\oddlayerDE}[6]{
  \coordinate (atop) at (#1,\frzTop+\frzExtra);
  \coordinate (abot) at (#1,\frzBot-\frzExtra);  
\foreach \ver/\pos/\val in {a1/5/{#2}, a3/3/{#3}, a5/1/{#4}, a7/0/{#5}, a9/-1/{#6}}
  \draw ({#1},\pos) node(\ver) {\val}; }
\newcommand{\evenlayerDE}[4]{
  \coordinate (atop) at (#1,\frzTop+\frzExtra);
  \coordinate (abot) at (#1,\frzBot-\frzExtra);  
  \foreach \ver/\pos/\val in {a2/4/{#2}, a4/2/{#3}, a6/0/{#4}}
    \draw ({#1},\pos) node(\ver) {\val}; }
\perEdge{0.5}; \perEdge{4.5};
\evenlayerDE{0}{$3$}{\vdots}{$n-1$}
\oddlayerDE{1}{$x$}{$2y$}{$mx$}{$x_{+}$}{$x_{-}$} \connect
\evenlayerDE{2}{$3$}{\vdots}{$n-1$} \connect 
\oddlayerDE{3}{$y$}{$2x$}{$my$}{$y_{+}$}{$y_{-}$} \connect 
\evenlayerDE{4}{$3$}{\vdots}{$n-1$} \connect 
\oddlayerDE{5}{$x$}{$2y$}{$mx$}{$x_{+}$}{$x_{-}$} \connect
\end{tikzpicture}
\] %\end{figure}
This also shows that the constant $\RRp$-frieze in \S\ref{sec:const} is unique, for $n$ even.

Now we consider friezes of type $D_{n}$ for $n$ odd, say $n=2m+1$.

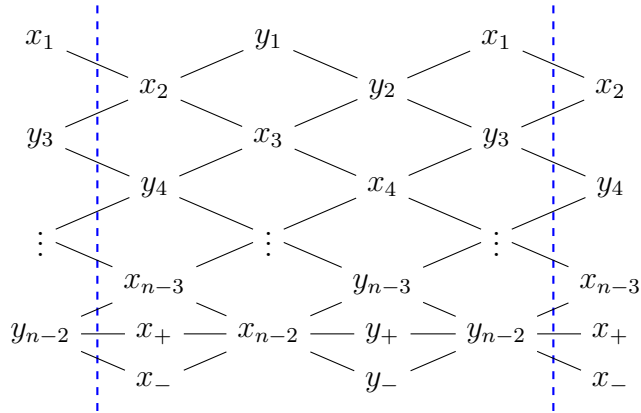
\[ %\begin{figure}[ht]
\begin{tikzpicture}[friezeM] \setfrzDO % ==  D Even  ==
\newcommand{\oddlayerDO}[6]{
  \coordinate (atop) at (#1,\frzTop+\frzExtra);
  \coordinate (abot) at (#1,\frzBot-\frzExtra);  
\foreach \ver/\pos/\val in {a1/5/{#2}, a3/3/{#3}, a5/1/{#4}, a7/0/{#5}, a9/-1/{#6}}
  \draw ({#1},\pos) node(\ver) {\val}; }
\newcommand{\evenlayerDO}[5]{
  \coordinate (atop) at (#1,\frzTop+\frzExtra);
  \coordinate (abot) at (#1,\frzBot-\frzExtra);  
  \foreach \ver/\pos/\val in {a0/6/{#2},a2/4/{#3}, a4/2/{#4}, a6/0/{#5}}
    \draw ({#1},\pos) node(\ver) {\val}; }
\perEdge{0.5}; \perEdge{4.5};
\evenlayerDO{0}{$x_1$}{$y_3$}{\vdots}{$y_{n-2}$}
\oddlayerDO{1}{$x_2$}{$y_4$}{$x_{n-3}$}{$x_{+}$}{$x_{-}$} \connect
\evenlayerDO{2}{$y_1$}{$x_3$}{\vdots}{$x_{n-2}$} \connect 
\oddlayerDO{3}{$y_2$}{$x_4$}{$y_{n-3}$}{$y_{+}$}{$y_{-}$} \connect 
\evenlayerDO{4}{$x_1$}{$y_3$}{\vdots}{$y_{n-2}$} \connect 
\oddlayerDO{5}{$x_2$}{$y_4$}{$x_{n-3}$}{$x_{+}$}{$x_{-}$} \connect
\end{tikzpicture}
\] %\end{figure}

The frieze equations are the same as in \eqref{eq:DEfrzeq},
but the parity difference means that we find that both even and odd levels are constant.
On the other hand, equations like \eqref{eq:c2k} apply to the odd levels and so
we find $c_{2k-1}=c_{2k+1}-2$, for $k=1,\ldots,m-1$,
but don't have a known initial value.

Then from the frieze equations, similarly to \eqref{eq:xyodd}, we get 
\begin{equation}\label{eq:coddeven}
  c_{2k}^2 = c_{2k-1}c_{2k+1}+1 = (c_{2k+1}-1)^2,
\end{equation}
so that $c_{2k}=c_{2k+1}-1$, since we are only looking for solutions in $\RRp$.
Thus we can start from $c_0=1$ to deduce that $c_k=k+1$ for all $k=1,\ldots,n-2$.
The only remaining equations are now
\[
x_{+}y_{+}  = x_{-}y_{-} = c_{n-2} + 1 = n,
\]
together with $x_{+}x_{-} = y_{+}y_{-}$, which then implies that $x_{-}=y_{+}$ and $y_{-}=x_{+}$.
Thus there is effectively one free parameter given by $\{(x,y)\in(\RRp)^2 : xy=n\}$ and the frieze is
\[ %\begin{figure}[ht]
\begin{tikzpicture}[friezeM] \setfrzDO % ==  D Even  ==
\newcommand{\oddlayerDO}[6]{
  \coordinate (atop) at (#1,\frzTop+\frzExtra);
  \coordinate (abot) at (#1,\frzBot-\frzExtra);  
\foreach \ver/\pos/\val in {a1/5/{#2}, a3/3/{#3}, a5/1/{#4}, a7/0/{#5}, a9/-1/{#6}}
  \draw ({#1},\pos) node(\ver) {\val}; }
\newcommand{\evenlayerDO}[5]{
  \coordinate (atop) at (#1,\frzTop+\frzExtra);
  \coordinate (abot) at (#1,\frzBot-\frzExtra);  
  \foreach \ver/\pos/\val in {a0/6/{#2},a2/4/{#3}, a4/2/{#4}, a6/0/{#5}}
    \draw ({#1},\pos) node(\ver) {\val}; }
\perEdge{0.5}; \perEdge{4.5};
\evenlayerDO{0}{$2$}{$4$}{\vdots}{$n-1$}
\oddlayerDO{1}{$3$}{$5$}{$n-2$}{$x$}{$y$} \connect
\evenlayerDO{2}{$2$}{$4$}{\vdots}{$n-1$} \connect 
\oddlayerDO{3}{$3$}{$5$}{$n-2$}{$y$}{$x$} \connect 
\evenlayerDO{4}{$2$}{$4$}{\vdots}{$n-1$} \connect 
\oddlayerDO{5}{$3$}{$5$}{$n-2$}{$x$}{$y$} \connect
\end{tikzpicture}
\] %\end{figure}
We get a constant frieze if and only if $x = y = \sqrt{n}$, so the constant frieze is unique. 

% ===============================================
\section{2-Periodic Friezes of Type $E$}
% ===============================================

Initial analysis, as in type $D_n$, shows that 2-periodic friezes of type $E_6$ and $E_8$ 
must actually be constant.
In the remaining case $E_7$, the same analysis shows which $\tau$-orbits must be constant,
so the general frieze will be as follows.

\begin{equation}\label{eq:E7frz}
\begin{tikzpicture}[friezeM, baseline=(bb.base)]
\pgfmathsetmacro{\frzTop}{3}
\pgfmathsetmacro{\frzBot}{-2}
\pgfmathsetmacro{\frzExtra}{0.7}
\newcommand{\connect}{
\foreach \x/\y in {a1/a2, a2/a3, a3/a4, a4/a5, a4/a7, a6/a7} \draw (\x)--(\y); }
\newcommand{\oddEdge}{
  \draw [edge] (atop)--(a1)--(a3)--(a5)--(a7)--(abot); }
\newcommand{\oddlayerES}[5]{
  \coordinate (atop) at (#1,\frzTop+\frzExtra);
  \coordinate (abot) at (#1,\frzBot-\frzExtra);  
  \foreach \ver/\pos/\val in {a1/3/{#2}, a3/1/{#3}, a5/0/{#4}, a7/-1/{#5}}
    \draw ({#1},\pos) node(\ver) {\val};
}
\newcommand{\evenlayerES}[4]{
  \coordinate (atop) at (#1,\frzTop+\frzExtra);
  \coordinate (abot) at (#1,\frzBot-\frzExtra);  
  \foreach \ver/\pos/\val in {a2/2/{#2}, a4/0/{#3}, a6/-2/{#4}}
    \draw ({#1},\pos) node(\ver) {\val}; 
}
\perEdge{0.5}; \perEdge{4.5};
\oddlayerES{0}{$x_1$}{$y_3$}{$y_7$}{$c_5$}
\evenlayerES{1}{$c_2$}{$c_4$}{$c_6$} \connect
\oddlayerES{2}{$y_1$}{$x_3$}{$x_7$}{$c_5$} \connect
\evenlayerES{3}{$c_2$}{$c_4$}{$c_6$} \connect
\oddlayerES{4}{$x_1$}{$y_3$}{$y_7$}{$c_5$} \connect
\evenlayerES{5}{$c_2$}{$c_4$}{$c_6$} \connect
\coordinate (bb) at (0,0);
\end{tikzpicture}
\end{equation}
The frieze equations are
\begin{equation}\label{eq:E7frzeq1}
\begin{aligned}
x_1y_1 &= c_2+ 1 \\
x_1y_3 = y_1x_3 &= c_2^2 - 1 \\
x_3y_3 &= c_2c_4 + 1 \\
x_7y_7 &= c_4+ 1 \\
x_3x_7 = y_3y_7 &= (c_4^2 -1)/c_5  \\
c_5^2 &= c_4c_6 + 1\\
c_6^2 &= c_5 + 1
\end{aligned}
\end{equation}
Setting $c_2=a+1$, we can solve the second (pair) to get $x_3=ax_1$ and $y_3=ay_1$,
so the first half of the equations reduce to
\begin{equation}\label{eq:E7frzeq2}
\begin{aligned}
x_1y_1 &= a+ 2 \\
c_4 &=  (x_3y_3-1)/c_2 = (a^2(a+2)-1)/(a+1) = a^2+a-1\\
x_7y_7 &= a(a+1) \\
x_1x_7 &= y_1y_7 \bigl(=\sqrt{a(a+1)(a+2)}\bigr)\\
\end{aligned}
\end{equation}
In particular, $a=x_3/x_1>0$ and $a^2+a-1=c_4>0$, so $a> \phi^{-1}$.
If the value of $a$ is determined, then there is one free parameter.
For the second half of the equations, we set $c_6=c$, giving
$c_5=c^2-1$ and $c_4=c(c^2-2)$,
and thus the constraint that $c>\sqrt{2}$.
Then, for consistency, we need that
\begin{equation}\label{eq:concy1}
\begin{aligned}
  a^2+a-1 &= c_4 = c(c^2-2), \\
  a\sqrt{a(a+1)(a+2)} &= x_3x_7= (c_4^2 -1)/c_5 = c^4-3c^2+1.
\end{aligned}
\end{equation}
Setting $x_1=y_1=\sqrt{a+2}$ will give a constant $\RRp$-frieze, 
for any value of $a$ satisfyiing the consistency conditions above.
Since we know from \S\ref{sec:const} that such a frieze is unique, we know that $a$ and $c$ must be given
by their values in that case, 
namely $a=2\phi$, $c=1+\phi$,
which one can check do satisfy \eqref{eq:concy1}.

Thus we end up with the following equations
\begin{equation}\label{eq:last-E7}
  x_1y_1 = 2\phi^2 = 2(1+\phi),  \quad
  x_7y_7 = 2\phi^4 = 2(2+3\phi), \quad
  x_1x_7 = y_1y_7 (= 2\phi^3)
\end{equation}
and we do indeed have one free parameter, similar to the case of $D_n$ ($n$ odd).
The remaining variables in \eqref{eq:E7frz} are given by
\begin{equation}\label{eq:last-E7-rest}
\begin{aligned}
  x_3 &= 2\phi x_1, \quad 
  y_3 = 2\phi y_1, \quad 
  c_6 = \phi^2 = 1+\phi, \quad 
  c_5 = 1+3\phi, \\
  c_2 &= \phi^3 = 1+2\phi, \quad 
  c_4 = 3\phi^3 = 3(1+2\phi).
\end{aligned}
\end{equation}

% ===============================================
\section{4-Periodic Friezes of Type $E_8$}
% ===============================================

As a final curiosity, we give the following example of a one-parameter family of 4-periodic $E_8$ friezes.
We can translate this frieze one step to get a second one-parameter family.
Translating two steps just gives $s\leftrightarrow t$.
 
\newcommand{\setfrzEA}{% == E 8 ==
\pgfmathsetmacro{\frzTop}{3}
\pgfmathsetmacro{\frzBot}{-4}
\pgfmathsetmacro{\frzExtra}{0.7}
\newcommand{\connect}{
\foreach \x/\y in {a1/a2, a1/a4, a3/a4, a4/a5, a5/a6, a6/a7, a7/a8} \draw (\x)--(\y); }
\newcommand{\oddEdge}{
  \draw [edge] (atop)--(a1)--(a3)--(a5)--(a7)--(abot); }
}
\newcommand{\oddlayerEA}[5]{
  \coordinate (atop) at (#1,\frzTop+\frzExtra);
  \coordinate (abot) at (#1,\frzBot-\frzExtra);  
  \foreach \ver/\pos/\val in {a1/2/{#2}, a3/1/{#3}, a5/-1/{#4}, a7/-3/{#5}}
    \draw ({#1},\pos) node(\ver) {\val};
}
\newcommand{\evenlayerEA}[5]{
  \coordinate (atop) at (#1,\frzTop+\frzExtra);
  \coordinate (abot) at (#1,\frzBot-\frzExtra);  
  \foreach \ver/\pos/\val in {a2/3/{#2}, a4/0/{#3}, a6/-2/{#4}, a8/-4/{#5}}
    \draw ({#1},\pos) node(\ver) {\val}; 
}

\begin{figure}[h]
\begin{tikzpicture}[xscale=2.0, yscale=1.0] 
  \setfrzEA % ==  E 8  ==
\evenlayerEA{-1}{$2+t$}{$3s+12+9t+2t^2$}{$3+3s+s^2$}{$1+s$}
\oddlayerEA{0}{$3+4t+t^2$}{$3+2t$}{$2t+7+4s+s^2$}{$2s+s^2$} \connect \oddEdge
\evenlayerEA{1}{$2+t$}{$3s+12+9t+2t^2$}{$3+3s+s^2$}{$1+s$} \connect
\oddlayerEA{2}{$5+2s+2t$}{$3+s+t$}{$9+3s+3t$}{$2+s+t$} \connect
\evenlayerEA{3}{$2+s$}{$3t+12+9s+2s^2$}{$3+3t+t^2$}{$1+t$} \connect
\oddlayerEA{4}{$3+4s+s^2$}{$3+2s$}{$2s+7+4t+t^2$}{$2t+t^2$} \connect \oddEdge
\evenlayerEA{5}{$2+s$}{$3t+12+9s+2s^2$}{$3+3t+t^2$}{$1+t$} \connect
\end{tikzpicture}
\caption{A 4-periodic frieze with extra reflection symmetry ($st=2$)}
\label{fig:oneparamE8}
\end{figure}
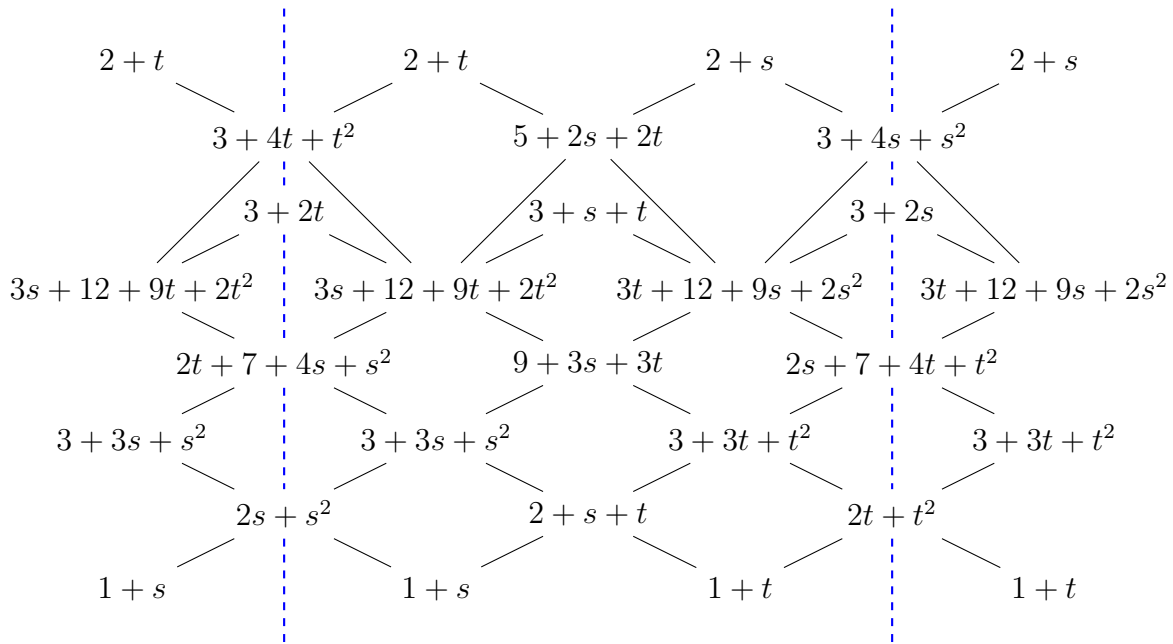

Setting $s=t=\sqrt{2}$ gives the constant frieze in \eqref{eq:E8const}.
Setting $(s,t)=(2,1)$ or $(1,2)$ in this or its translate gives 
the four well-known integer friezes of type $E_8$ with no 1's
(see e.g.~\cite[Example~5.15]{BFGST}).
We do not know an explicit classification of the 4-periodic $E_8$ friezes,
but one might conjecture that there should be a two-parameter family
containing the two one-parameter families given here.
Note that this is consistent with \cite[Thm~3.12]{SGL}, 
in that the derivative $D\tau$ of $\tau$ at the fixed point has eigenvalues $i$ and $-i$
with multiplicity one and hence there is a two dimensional subspace on which
$D\tau$ has order 4.

\goodbreak
% ===============================================

% ===============================================
\end{document}